\documentclass[twocolumn]{autart}    

\usepackage{graphics} 
\usepackage{epsfig} 
\usepackage{mathptmx} 
\usepackage{times} 
\usepackage{amsmath} 
\usepackage{amssymb}  
\usepackage{epstopdf}
\usepackage{subfigure}
\usepackage{amssymb}
\usepackage{color}

\graphicspath{{./Figures/}}
\allowdisplaybreaks[4]

\newcommand{\norme}[1]{\left\Vert #1\right\Vert}

\begin{document}

\begin{frontmatter}

\title{Boundary Feedback Stabilization of a Flexible Wing Model under Unsteady Aerodynamic Loads\thanksref{footnoteinfo}} 

\thanks[footnoteinfo]{This paper was not presented at any IFAC
meeting. Corresponding author G.~Zhu. }

\author[Polymtl]{Hugo Lhachemi}\ead{hugo.lhachemi@polymtl.ca},        
\author[Polymtl]{David Saussi\'{e}}\ead{d.saussie@polymtl.ca},        
\author[Polymtl]{Guchuan Zhu}\ead{guchuan.zhu@polymtl.ca}             

\address[Polymtl]{\'{E}cole Polytechnique de Montr\'{e}al, Montr\'{e}al, Canada}  


\begin{abstract}                          
This paper addresses the boundary stabilization of a flexible wing model, both in bending and twisting displacements, under unsteady aerodynamic loads, and in presence of a store. The wing dynamics is captured by a distributed parameter system as a coupled Euler-Bernoulli and Timoshenko beam model. The problem is tackled in the framework of semigroup theory, and a Lyapunov-based stability analysis is carried out to assess that the system energy, as well as the bending and twisting displacements, decay exponentially to zero. The effectiveness of the proposed boundary control scheme is evaluated based on simulations.
\end{abstract}

\end{frontmatter}

\section{Introduction}
Modern aerospace systems such as aircraft, unmanned aerial vehicles (UAVs), and microaerial vehicles are subject to stringent performance requirements including high maneuverability and extended autonomy. The global trend to achieve the required level of performance consists in reducing the mass of the system by a massive integration of composite materials. However, it results in a decrease of the structure rigidity. In particular, lightweight flexible wings are subject to stronger aeroelastic phenomena which are the result from interactions between aerodynamic, elastic and inertial forces. Such phenomena can significantly degrade the performance of an aircraft by introducing undesired couplings between the flexible modes and the flight dynamics~\cite{Shearer2007,Su2010}, and may also jeopardize the integrity of its structure~\cite{Mukhopadhyay2003}. These phenomena can be amplified in the case of a store located under the wing with the emergence of the so-called store-induced oscillations~\cite{beran2004studies,Bialy2016}. Therefore, the active control of aeroelastic phenomena has become a topic of primary interest.

One of the most noticeable contributions for the control of aeroelastic phenomena is the Benchmark Active Control Technology (BACT) wind-tunnel model developed by NASA Langley Research Center~\cite{Scott2000}. The BACT is modeled as a two-degree-of-freedom aeroelastic wing section capturing the first bending and twisting modes of a flexible wing. The control design strategy of the BACT for flutter suppression, including experimental tests, has been widely investigated in the literature~\cite{Bhoir2003,Ko1999,Mukhopadhyay2000Transonic}. Nevertheless, the BACT cannot fully represent the dynamics of real flexible wings. Indeed, the flexible wing can be more accurately modeled by a distributed parameter system of two coupled partial differential equations (PDEs) describing the dynamics in bending and twisting displacements respectively~\cite{Bialy2016,Zhang2005,Ziabar2010}.

The study on flexible structures described by distributed systems and their interactions with the flow-field has attracted many attention in the last decades~\cite{stanewsky2001adaptive}. The bending dynamics of a panel evolving in different flow-field regimes have been studied for clamped~\cite{chueshov2012generation,lasiecka2016feedback} and clamped-free~\cite{chueshov2016nonlinear} boundary conditions in case of a distributed velocity feedback. The coupled Euler-Bernoulli and Timoshenko beam model, describing both undamped bending and torsion flexible displacements, has also been investigated for self-straining actuators employed as boundary control inputs without~\cite{balakrishnan2004spectral} and with~\cite{balakrishnan2001subsonic,balakrishnan2003toward} an external load generated by the flow-field.

This paper addresses the boundary stabilization problem of a flexible wing whose dynamics are captured by a coupled Euler-Bernoulli and Timoshenko beam model in the presence of a store located at the wing tip. Unlike the self-straining actuation setup considered in~\cite{balakrishnan2001subsonic,balakrishnan2003toward,balakrishnan2004spectral}, the actuation scheme consists in flaps located at the wing tip to locally generate lift force and torsional momentum, resulting in distinct boundary conditions of the coupled PDEs. Furthermore, the model considered in this work includes the contribution of the Kelvin-Voigt damping~\cite{zhang2011spectrum} in both bending and twisting axes. This model is a linear and damped version of the one presented in~\cite{Bialy2016}, while the aerodynamic loads are supposed to be unsteady. A similar problem, namely a flapping wing UAV, is considered in~\cite{Paranjape2013}. The model used includes the contribution of the Kelvin-Voigt damping, while assuming that aerodynamic loads are unknown but bounded. The method of backstepping is used for the boundary control of the spatial integral of the state variables to track the net aerodynamic forces on the wing. The same model is considered in~\cite{He2017} for which a Lyapunov-based stabilization control is developed to achieve bounded bending and twisting deflections in the presence of aerodynamic load disturbances. It is worth noting that as pointed out in \cite{Curtain2009}, an Euler-Bernoulli beam model with Kelvin-Voigt damping may not be well-posed if the boundary conditions do not explicitly include the Kelvin-Voigt damping term in a correct manner. Therefore, although the existence of Kelvin-Voigt damping may intuitively be helpful for system stabilization, a rigorous well-posedness analysis is still needed to guarantee the expected behavior of the considered system which remains a more complex setting than a single beam. This constitutes one of the main motivations of the present work.

It should be noticed that a commonly used assumption for Lyapunov-based designs in the works~\cite{Bialy2016,He2017,Paranjape2013} is that either the system energy or the aerodynamic loads should be bounded. Furthermore it is also assumed the existence and the regularity of the system trajectories and their partial derivatives up to a certain order. These assumptions, which can be justified by physical intuitions~\cite{Queiroz2012,Queiroz2002}, can considerably simplify closed-loop stability analysis. However, they imply the well-posedness of the underlying PDEs, which is a quite strong condition. The main objective of this work is to show that these assumptions can be relaxed. To this aim, we formulate the problem under an abstract form that allows the application of the semigroup theory~\cite{Curtain2012,Pazy2012}. Due to the presence of the store, the boundary conditions related to the control strategy take the form of ODEs \cite{Guo2002}. It results in a system in abstract form composed of two coupled PDEs and two coupled ODEs. We show that the closed-loop system with the proposed boundary control admits a $C_{0}$-semigroup and is well-posed. The closed-loop stability is derived from a Lyapunov-based analysis, which shows that the above $C_{0}$-semigroup is exponentially stable. The results of this work allow confirming the validity of most existing control schemes reported in the literature for similar settings under even much less restrictive conditions.

The remainder of the paper is organized as follows. The wing model, along with the associated abstract form, are introduced in Section~\ref{sec: model}. The well-posedness of the problem is analyzed in Section~\ref{sec: well-posedness} in the framework of semigroup theory. Then, a Lyapunov-based analysis is carried out in Section~\ref{sec: stability} to assess that the system energy, as well as bending and twisting displacements, exponentially decay to zero. Finally, numerical simulations are presented in Section~\ref{sec: simulations} to illustrate the performance of the closed-loop system.

\textbf{Notations \cite{Leoni2009,Royden1988}:} $\mathbb{R}_+$ and $\mathbb{R}_+^*$ denote the sets of non-negative and positive real numbers, respectively. Let $L^2(0,l)$ be the set of Lebesgue squared integrable real-valued functions over $(0,l)$ endowed with its natural norm denoted by $\norme{\cdot}_{L^2(0,l)}$. For any $m\in\mathbb{N}$, $H^m(0,l)$ denotes the usual Sobolev space, which is defined as the set of $f\in L^2(0,l)$, such that $f$ admits $m$ successive weak derivatives, denoted by $f',f'',\ldots,f^{(m)}$, in $L^2(0,l)$. Denoting by $\mathrm{AC}[0,l]$ the set of all absolutely continuous functions on $[0,l]$, $H^1(0,l) \subset \mathrm{AC}[0,l]$ in the sense that for any $f\in H^1(0,l)$, there exists a unique absolutely continuous function $g\in\mathrm{AC}[0,l]$ such that $f=g$ in $H^1(0,l)$. We note $H_L^m(0,l) = \{f \in H^m(0,l) \; : \; f(0)=f'(0)=\ldots=f^{(m-1)}(0)=0\}$. For a given normed vector spaces $(E,\norme{\cdot}_{E})$, $\mathcal{L}(E)$ denotes the space of bounded linear transformations from $E$ to $E$. The range of a given operator $\mathcal{A}$ is denoted by $R(\mathcal{A})$ while its resolvent set is denoted by $\rho(\mathcal{A})$. The successive partial derivatives of a sufficiently regular function $f$ are denoted in subscript, e.g., $f_{ty}$ stands for $\partial^2 f / (\partial t \partial y)$.

\section{Problem Setting and Boundary Control Law}\label{sec: model}
\subsection{Flexible wing model}
Let $l\in\mathbb{R}_+^*$ be the length of the wing, $\rho\in\mathbb{R}_+^*$ the mass per unit of span, $I_w\in\mathbb{R}_+^*$ the moment of inertia per unit length, $EI\in\mathbb{R}_+^*$ (resp. $GJ\in\mathbb{R}_+^*$) the bending (resp. torsional) stiffness, $\eta_\omega\in\mathbb{R}_+^*$ (resp. $\eta_\phi\in\mathbb{R}_+^*$) the bending (resp. torsional) Kelvin-Voigt damping coefficient, and $x_c \in \mathbb{R}$ the distance between the wing center of gravity and the elastic axis of the wing. The store at the wing tip is characterized by its mass $m_s\in\mathbb{R}_+^*$ and its moment of inertia $J_s\in\mathbb{R}_+^*$. We define the two following symmetric definite positive matrices:
\begin{equation}
M \triangleq
\begin{bmatrix}
\rho & \rho x_c \\ \rho x_c & I_w^*
\end{bmatrix}
, \qquad
M_s \triangleq
\begin{bmatrix}
m_s & m_s x_c \\ m_s x_c & J_s^*
\end{bmatrix} ,
\end{equation}
with $I_w^* \triangleq I_w + \rho x_c^2$ and $J_s^* \triangleq J_s + m_s x_c^2$. Introducing $c_\omega = \sqrt{EI/\rho}$ and $c_\phi = \sqrt{GJ/I_w}$, the bending and twisting dynamics are described by the following coupled PDEs \cite{Bialy2016,He2017,Paranjape2013}:
\begin{equation} \label{eq: studied EDP}
M
\begin{bmatrix}
\omega_{tt} \\ \phi_{tt}
\end{bmatrix}
+
\begin{bmatrix}
\rho c_\omega^2 (\omega_{yy} + \eta_\omega \omega_{tyy})_{yy} \\ - I_w c_\phi^2 (\phi_{y} + \eta_\phi \phi_{ty})_{y}
\end{bmatrix}
=
\begin{bmatrix}
F_a \\ M_a
\end{bmatrix} ,
\end{equation}
where the functions $\omega : [0,l] \times \mathbb{R}_+ \rightarrow \mathbb{R}$ and $\phi : [0,l] \times \mathbb{R}_+ \rightarrow \mathbb{R}$ denote, respectively, the bending and twisting displacements at the location $y\in[0,l]$ along the wing span and at time $t \geq 0$ and $F_a : [0,l] \times \mathbb{R}_+ \rightarrow \mathbb{R}$ and $M_a : [0,l] \times \mathbb{R}_+ \rightarrow \mathbb{R}$ denote, respectively, the aerodynamic lift force and pitching moment applied at the location $y\in[0,l]$ and at time $t \geq 0$. They are expressed under the following unsteady form:
\begin{equation}\label{eq: aero efforts}
\begin{bmatrix}
F_a \\ M_a
\end{bmatrix}
\triangleq
\begin{bmatrix}
\alpha_\omega \phi + \beta_\omega \phi_t + \gamma_\omega \omega_t \\
\alpha_\phi \phi + \beta_\phi \phi_t + \gamma_\phi \omega_t 
\end{bmatrix} ,
\end{equation}
where $\alpha_\omega,\beta_\omega,\gamma_\omega,\alpha_\phi,\beta_\phi,\gamma_\phi\in\mathbb{R}_+$. This model, commonly employed in finite dimension~\cite{Bhoir2003,Ko1999,Mukhopadhyay2000Transonic}, is a trade-off between the used steady form of~\cite{Bialy2016} and the unmodeled black-box representation in~\cite{He2017,Paranjape2013}. The boundary conditions for the tip-based control scheme, in the presence of a store~\cite{Bialy2016}, considered in this work are such that, for any $t\geq0$,
\begin{subequations}
\begin{eqnarray}
\omega(0,t) & = & \omega_{y}(0,t)=\omega_{yy}(l,t) = \phi(0,t) = 0 , \label{eq: EDP boundary zero} \\
M_s \begin{bmatrix} \omega_{tt}(l,t) \\ \phi_{tt}(l,t) \end{bmatrix} & = & \begin{bmatrix} L_\mathrm{tip}(t) + \rho c_\omega^2 (\omega_{yy} + \eta_\omega \omega_{tyy})_{y}(l,t) \\ M_\mathrm{tip}(t) - I_w c_\phi^2 (\phi_{y} + \eta_\phi \phi_{ty})(l,t) \end{bmatrix} , \label{eq: EDP boundary flexion and torsion}
\end{eqnarray}
\end{subequations}
where $L_\mathrm{tip}:\mathbb{R}_+\rightarrow\mathbb{R}$ and $M_\mathrm{tip}:\mathbb{R}_+\rightarrow\mathbb{R}$ are the tip control inputs. More precisely, $L_\mathrm{tip}(t)$ and $M_\mathrm{tip}(t)$ denote the aerodynamic lift force and pitching moment generated at time $t$ by the flaps located at the wing tip. Finally, the initial conditions are given, for any $y\in(0,l)$, by
\begin{subequations}
\begin{align}
\omega(y,0) & = \omega_0(y) ,\;\; \omega_t(y,0) = \omega_{t0}(y) , \\
\phi(y,0) & = \phi_0(y) ,\;\; \phi_t(y,0) = \phi_{t0}(y) .
\end{align}
\end{subequations}

\subsection{Boundary control law}
Introducing the system energy defined by
\begin{align}
E \triangleq \dfrac{1}{2} & \int_{0}^{l} \rho c_\omega^2 \omega_{yy}^2 + I_w c_\phi^2 \phi_{y}^2
+ \begin{bmatrix}
\omega_{t} \\ \phi_{t}
\end{bmatrix}^\top
M
\begin{bmatrix}
\omega_{t} \\ \phi_{t}
\end{bmatrix}
\mathrm{d}y \label{eq: definition energy} \\
& + \dfrac{1}{2}
\begin{bmatrix}
\omega_t(l,\cdot) \\ \phi_t(l,\cdot)
\end{bmatrix}^\top
M_s
\begin{bmatrix}
\omega_t(l,\cdot) \\ \phi_t(l,\cdot)
\end{bmatrix} \nonumber ,
\end{align}
the control problem investigated in this paper is formalized as follows.
\begin{prob}\label{eq: boundary control pb}
The boundary control objective is twofold.
\begin{enumerate}
\item To guarantee that the system energy exponentially decays to zero, i.e., there exist $K_E,\Lambda\in\mathbb{R}_+^*$ such that
\begin{equation*}
\forall t \geq 0,  \;
E(t) \leq K_E E(0) \exp(-\Lambda t) .
\end{equation*}
\item To guarantee that both bending and twisting displacements converge exponentially and uniformly over the wing span to zero.
\end{enumerate}
\end{prob}
In particular, it will be shown that if the first objective in Problem~\ref{eq: boundary control pb} is satisfied, the second one follows in the sense that there exist $K_\omega,K_\phi\in\mathbb{R}_{+}^{*}$ such that for any $t \geq 0$,
\begin{equation*}
\begin{split}
\norme{\omega(\cdot,t)}_\infty & \leq K_\omega \sqrt{E(0)} \exp(-\Lambda t /2) , \\
\norme{\phi(\cdot,t)}_\infty & \leq K_\phi \sqrt{E(0)} \exp(-\Lambda t /2) , \\
\end{split}
\end{equation*}
where $\norme{\cdot}_\infty$ is the uniform norm for real-valued functions defined over $[0,l]$, i.e., $\norme{f}_\infty \triangleq \sup\left\{ \lvert f(y) \rvert : y\in[0,l] \right\}$.

In control design and implementation, we make the following assumption.
\begin{assum}
It is assumed that $\omega(l,\cdot)$, $\omega_{t}(l,\cdot)$, $\phi(l,\cdot)$, and $\phi_{t}(l,\cdot)$ are measured and available for feedback control.
\end{assum}
Note that the above assumptions are commonly used in the existing literature. In practice, the point-wise displacements of the structure at the wing tip can be measured by piezoelectric bending and torsion transducers, while their time derivative can be obtained by numerical methods.

The proposed boundary stabilizing control is formed by two proportional-derivative (PD) controllers:
\begin{subequations}\label{eq: boundary control}
\begin{align}
L_\mathrm{tip}(t) & = - k_1 \left[ \omega_{t}(l,t) + \epsilon_1 \omega(l,t) \right] , \label{eq: bf L_tip} \\
M_\mathrm{tip}(t) & = - k_2 \left[\phi_{t}(l,t) + \epsilon_2 \phi(l,t) \right] , \label{eq: bf M_tip}
\end{align}
\end{subequations}
for any $t\geq0$, where $k_1,k_2 \in\mathbb{R}_+^*$ are tunable controller gains and $\epsilon_1,\epsilon_2 \in\mathbb{R}_+^*$ are two parameters to be determined. In the remainder of this paper, we show that this boundary control solves Problem~\ref{eq: boundary control pb}.

\subsection{Closed-loop system in abstract form}
In order to study the properties of the closed-loop system, the problem is rewritten in abstract form. In this context, we introduce the following real Hilbert space:
\begin{equation*}
\begin{split}
\mathcal{H} = H_L^2(0,l) \times L^2(0,l) \times H_L^1(0,l) \times L^2(0,l) \times \mathbb{R} \times \mathbb{R} ,
\end{split}
\end{equation*}
endowed with the inner product $\left\langle \cdot,\cdot \right\rangle_{\mathcal{H},1}$ defined for any $X_i=(f_i,g_i,h_i,z_i,\zeta_{\omega,i},\zeta_{\phi,i})\in\mathcal{H}$, $i\in\{1,2\}$, by
\begin{equation*}
\begin{split}
& \langle X_1,X_2 \rangle_{\mathcal{H},1} \\
\triangleq & \int_0^l [ \rho c_\omega^2 f_1''(y) f_2''(y) + I_w c_\phi^2 h_1'(y) h_2'(y) ] \mathrm{d}y \\
& \phantom{\triangleq\;} + \int_0^l
\begin{bmatrix}
g_1(y) \\ z_1(y)
\end{bmatrix}^\top
M
\begin{bmatrix}
g_2(y) \\ z_2(y)
\end{bmatrix}
\mathrm{d}y
+ \begin{bmatrix} \zeta_{\omega,1} \\ \zeta_{\phi,1} \end{bmatrix}^\top M_s \begin{bmatrix} \zeta_{\omega,2} \\ \zeta_{\phi,2} \end{bmatrix} .
\end{split}
\end{equation*}
The induced norm $\norme{\cdot}_{\mathcal{H},1}$ is such that the energy of the wing defined by (\ref{eq: definition energy}) can be expressed for any $t \geq 0$ as
\begin{equation*}
E(t) = \dfrac{1}{2} \norme{(\omega(\cdot,t) , \omega_t(\cdot,t) , \phi(\cdot,t) , \phi_t(\cdot,t) , \omega(l,t) , \phi(l,t))}_{\mathcal{H},1}^2 .
\end{equation*}

In view of equations (\ref{eq: studied EDP}), the boundary conditions (\ref{eq: EDP boundary zero}-\ref{eq: EDP boundary flexion and torsion}), and the boundary control (\ref{eq: bf L_tip}-\ref{eq: bf M_tip}), we introduce the operator $\mathcal{A}$ defined on
{\allowdisplaybreaks[4]
\begin{align}
D(\mathcal{A}) \triangleq \{ & (f,g,h,z,\zeta_\omega,\zeta_\phi)\in\mathcal{H} \, : \label{eq: def domain A} \\
& \, g \in H_L^2(0,l) ,\, z \in H_L^1(0,l) ,\, \nonumber \\
& f''+\eta_\omega g'' \in H^2(0,l) ,\, h'+\eta_\phi z' \in H^1(0,l) ,\, \nonumber \\
& (f''+\eta_\omega g'')(l)=0 , \zeta_\omega = g(l) , \, \zeta_\phi = z(l) ,\, \nonumber \\
& f,f',g,g',h,z,(f''+\eta_\omega g'')\in\mathrm{AC}[0,l] ,\, \nonumber \\
& (f''+\eta_\omega g'')',(h'+\eta_\phi z')\in\mathrm{AC}[0,l] \} , \nonumber
\end{align}
}
by $\mathcal{A} = \mathcal{A}_1 + \mathcal{A}_2$ with $D(\mathcal{A}_1)=D(\mathcal{A})$, $D(\mathcal{A}_2) = \mathcal{H}$, $\mathcal{A}_1(f,g,h,z,\zeta_\omega,\zeta_\phi) = (g,\tilde{g},z,\tilde{z},\tilde{\zeta}_\omega,\tilde{\zeta}_\phi)$ where
\begin{subequations}
\begin{align}
\begin{bmatrix} \tilde{g} \\ \tilde{z} \end{bmatrix}
& \triangleq - M^{-1} \begin{bmatrix} \rho c_\omega^2(f''+\eta_\omega g'')'' \\ - I_w c_\phi^2 (h' + \eta_\phi z')' \end{bmatrix} , \label{eq: def A1 - part1} \\
\begin{bmatrix} \tilde{\zeta}_\omega \\ \tilde{\zeta}_\phi \end{bmatrix}
& \triangleq M_s^{-1} \begin{bmatrix} \rho c_\omega^2 (f''+\eta_\omega g'')'(l) - k_1 (g(l)+\epsilon_1 f(l)) \\ - I_w c_\phi^2 (h'+\eta_\phi z')(l) - k_2 (z(l)+\epsilon_2 h(l)) \end{bmatrix} , \label{eq: def A1 - part2}
\end{align}
\end{subequations}
and $\mathcal{A}_2(f,g,h,z,\zeta_\omega,\zeta_\phi) = (0,\tilde{F}_a,0,\tilde{M}_a,0,0)$ where
\begin{equation*}
\begin{bmatrix}
\tilde{F}_a \\ \tilde{M}_a
\end{bmatrix}
\triangleq
M^{-1}
\begin{bmatrix}
F_a \\ M_a
\end{bmatrix}
=
M^{-1}
\begin{bmatrix}
\alpha_\omega h + \beta_\omega z + \gamma_\omega g \\
\alpha_\phi h + \beta_\phi z + \gamma_\phi g ,
\end{bmatrix}.
\end{equation*}
Then, the evolution equation in abstract form is given by
\begin{equation}\label{eq: cauchy problem}
\left\{\begin{split}
\dfrac{\mathrm{d} X}{\mathrm{d} t}(t) & = \mathcal{A} X(t) ,\; t>0  ,\\
X(0) & = X_0 \in D(\mathcal{A}) ,
\end{split}\right.
\end{equation}
where $X(t) = \left( \omega(\cdot,t) , \omega_t(\cdot,t) , \phi(\cdot,t) , \phi_t(\cdot,t) , \omega_t(l,t) , \phi_t(l,t) \right)$ is the state vector and $X_0 = \left( \omega_0 , \omega_{t0} , \phi_0 , \phi_{t0} , \omega_{t0}(l) , \phi_{t0}(l) \right)$ is the initial condition.
\begin{rem}
The boundary condition $\omega_{yy}(l,t)=0$ (see (\ref{eq: EDP boundary zero})) implies $\omega_{tyy}(l,t)=0$, which provides $\omega_{yy}(l,t)+\eta_\omega \omega_{tyy}(l,t)=0$. Conversely, $\omega_{yy}(l,t)+\eta_\omega \omega_{tyy}(l,t)=0$ with the initial condition $\omega_{yy}(l,0)=0$ implies $\omega_{yy}(l,t)=0$ for any $t\geq0$. It motivates the introduction of the boundary constraint $(f''+\eta_\omega g'')(l)=0$ in (\ref{eq: def domain A}).
\end{rem}
The proof of the closed-loop exponential stability consists in two main steps:
\begin{enumerate}
\item to show that $\mathcal{A}$ is the infinitesimal generator of a $C_0$-semigroup $T(t)$ on $\mathcal{H}$;
\item to show that the $C_0$-semigroup $T(t)$ is exponentially stable.
\end{enumerate}
With this approach, the regularity properties of the closed-loop trajectory are deduced from the well-posedness assessment, and the exponential energy decay of the system in closed loop is confirmed by Lyapunov stability analysis. The details of the proof are presented in the two next sections.

\section{Well-posedness}\label{sec: well-posedness}
To assess that the Cauchy problem (\ref{eq: cauchy problem}) is well-posed, it is necessary to study the properties of the operator $\mathcal{A}$~\cite[Chap.4, Th.1.3.]{Pazy2012}. In the upcoming developments, the following versions of the Poincar\'{e}'s and Agmon's inequalities will be used.
\begin{lem}\label{lemma: Poincare's and Agmon's inequalities}~\cite{Hardy1952,Krstic2008}
For any $f\in H^1(0,l)$ such that $f(0)=0$, the Poincar\'{e}'s inequality ensures that
\begin{equation*}
\norme{f}_{L^2(0,l)}^2 \leq \dfrac{4 l^2}{\pi^2} \norme{f'}_{L^2(0,l)}^2 ,
\end{equation*}
while the Agmon's inequality provides
\begin{equation*}
\norme{f}_\infty^2 \leq 2 \norme{f}_{L^2(0,l)} \norme{f'}_{L^2(0,l)} .
\end{equation*}
\end{lem}
\subsection{Introduction of a second inner product on $\mathcal{H}$}
In order to study both the well-posedness and the stability properties of the abstract Cauchy problem (\ref{eq: cauchy problem}), it will be useful to consider a second inner product on $\mathcal{H}$. Such an approach is generally employed to enforce a dissipative property of the studied operator in an adequate Hilbert space and belongs to the framework of energy multipliers (see, e.g., \cite{krstic2008output,Luo2012}). Let $\epsilon_1,\epsilon_2\in\mathbb{R}_+^*$ be the constants involved in the control law (\ref{eq: bf L_tip}-\ref{eq: bf M_tip}) and $\Psi : \mathcal{H} \times \mathcal{H} \rightarrow \mathbb{R}$ be defined for any $X_i=(f_i,g_i,h_i,z_i,\zeta_{\omega,i},\zeta_{\phi,i})\in\mathcal{H}$, $i\in\{1,2\}$, by
\begin{align*}
\Psi(X_1,X_2) &
= \int_0^l
\begin{bmatrix}
g_1 \\ z_1
\end{bmatrix}^\top
M
\begin{bmatrix}
\epsilon_1 f_2 \\ \epsilon_2 h_2
\end{bmatrix}
+
\begin{bmatrix}
g_2 \\ z_2
\end{bmatrix}^\top
M
\begin{bmatrix}
\epsilon_1 f_1 \\ \epsilon_2 h_1
\end{bmatrix} \mathrm{d}y \\
& \phantom{=\;} + \begin{bmatrix}
\zeta_{\omega,1} \\ \zeta_{\phi,1}
\end{bmatrix}^\top
M_s
\begin{bmatrix}
\epsilon_1 f_2(l) \\ \epsilon_2 h_2(l)
\end{bmatrix}
+
\begin{bmatrix}
\zeta_{\omega,2} \\ \zeta_{\phi,2}
\end{bmatrix}^\top
M_s
\begin{bmatrix}
\epsilon_1 f_1(l) \\ \epsilon_2 h_1(l)
\end{bmatrix} .
\end{align*}
Then, let $\left\langle\cdot,\cdot\right\rangle_{\mathcal{H},2}:\mathcal{H}\times\mathcal{H} \rightarrow \mathbb{R}$ be defined by $\langle X_1 , X_2 \rangle_{\mathcal{H},2} = \langle X_1 , X_2 \rangle_{\mathcal{H},1} + \Psi(X_1,X_2)$. Finally, let $K_{m,1},K_{m,2}\in\mathbb{R}_+^*$ be defined by
\begin{align*}
K_{m,1} = &
\min \left( \dfrac{\pi^4 \rho c_\omega^2}{4 l^3 (1+|x_c|) \{ 4l \rho^{3/2} c_\omega + \pi^2 m_s \} } ,
\dfrac{c_\omega \lambda_m(M)}{2\sqrt{\rho}} , \right. \\
& \left. \phantom{===========} \dfrac{c_\omega \lambda_m(M)}{2\sqrt{\rho}|x_c|} , \dfrac{\lambda_m (M_s)}{2 m_s} , \dfrac{\lambda_m (M_s)}{2 m_s |x_c|} \right) , \\
K_{m,2} = &
\min \left( \dfrac{\pi^2 I_w c_\phi^2}{4 l^2 ( I_w^* + \rho |x_c|) \sqrt{I_w} c_\phi + \pi^2 l (J_s^*+m_s|x_c|)} , \right. \\
& \phantom{=} \left. \dfrac{\sqrt{I_w} c_\phi \lambda_m(M)}{2 \rho |x_c|} , \dfrac{\sqrt{I_w} c_\phi \lambda_m(M)}{2 I_w^*} , \dfrac{\lambda_m(M_s)}{2m_s|x_c|} , \dfrac{\lambda_m(M_s)}{2 J_s^*} \right) ,
\end{align*}
where $\lambda_m(\cdot)$ denotes the smallest eigenvalue.
\begin{lem}\label{lemma: inner product 2}
For any given $\epsilon_i \in (0,K_{m,i})$ with $i \in \{1,2\}$, $\left\langle\cdot,\cdot\right\rangle_{\mathcal{H},2}$ is an inner product for $\mathcal{H}$. Furthermore, the induced norm, denoted by $\norme{\cdot}_{\mathcal{H},2}$, is equivalent to $\norme{\cdot}_{\mathcal{H},1}$.
\end{lem}
\textbf{Proof.} Let $\epsilon_i = \alpha_i K_{m,i}$ with $\alpha_i \in (0,1)$ for $i \in \{1,2\}$. Note first that $\left\langle\cdot,\cdot\right\rangle_{\mathcal{H},2}$ is bilinear and symmetric. For any $X=(f,g,h,z,\zeta_\omega,\zeta_\phi)\in\mathcal{H}$, applying Young's\footnote{$\forall a,b\in\mathbb{R}_+$, $\forall r\in\mathbb{R}_+^*$, $ab \leq a^2/(2r) + r b^2/2$.}, Schwartz's, and Poincar\'{e}'s inequalities, it yields
\begin{align*}
& \vert \Psi(X,X) \vert \\
\leq & \epsilon_1 \dfrac{4l^3}{\pi^2} (1+|x_c|) \left\{ \dfrac{4l \sqrt{\rho}}{\pi^2 c_\omega} + \dfrac{m_s}{\rho c_\omega^2} \right\} \int_0^l \rho c_\omega^2 f''(y)^2 \mathrm{d}y \\
& + \epsilon_2 l \left\{ \dfrac{4l}{\pi^2} \dfrac{I_w^* + \rho |x_c|}{\sqrt{I_w} c_\phi} + \dfrac{J_s^* + m_s |x_c|}{I_w c_\phi^2} \right\} \int_0^l I_w c_\phi^2 h'(y)^2 \mathrm{d}y \\
& + \left( \epsilon_1 \dfrac{\sqrt{\rho}}{c_\omega} + \epsilon_2 \dfrac{\rho |x_c|}{\sqrt{I_w}c_\phi} \right) \int_0^l g(y)^2 \mathrm{dy} \\
& + \left( \epsilon_1 \dfrac{\sqrt{\rho} |x_c|}{c_\omega} + \epsilon_2 \dfrac{I_w^*}{\sqrt{I_w}c_\phi} \right) \int_0^l z(y)^2 \mathrm{dy} \\
& + (\epsilon_1 m_s + \epsilon_2 m_s |x_c|) \zeta_\omega^2 + (\epsilon_1 m_s |x_c| + \epsilon_2 J_s^*) \zeta_\phi^2 .
\end{align*}
Recalling that for any symmetric matrix $S\in\mathbb{R}^{n \times n}$ and any vector $x\in\mathbb{R}^n$, $x^\top S x \geq \lambda_m(S) x^\top x$, it provides $\vert \psi(X,X) \vert \leq \alpha \norme{X}_{\mathcal{H},1}^2$ with $\alpha = \max(\alpha_1,\alpha_2) \in (0,1)$. From the definition of $\langle\cdot,\cdot\rangle_{\mathcal{H},2}$, we deduce that:
\begin{equation}\label{eq: norm equivalent}
( 1 - \alpha ) \norme{X}_{\mathcal{H},1}^2
\leq
\left\langle X , X \right\rangle_{\mathcal{H},2}
\leq
( 1 + \alpha ) \norme{X}_{\mathcal{H},1}^2 .
\end{equation}
Then, $\left\langle\cdot,\cdot\right\rangle_{\mathcal{H},2}$ is positive definite and hence, it defines an inner product on $\mathcal{H}$. Denoting by $\norme{\cdot}_{\mathcal{H},2}$ the induced norm, (\ref{eq: norm equivalent}) shows that $\norme{\cdot}_{\mathcal{H},1}$ and $\norme{\cdot}_{\mathcal{H},2}$ are equivalent. \hfill $\qed$

From Lemma~\ref{lemma: inner product 2} the following corollary holds.
\begin{cor}
$(\mathcal{H},\left\langle\cdot,\cdot\right\rangle_{\mathcal{H},2})$ is a real Hilbert space.
\end{cor}
\begin{rem}
The utility of the second inner product $\left\langle\cdot,\cdot\right\rangle_{\mathcal{H},2}$ for assessing the well-posedness of the abstract Cauchy problem (\ref{eq: cauchy problem}) will appear clearly in the proof of Lemma~\ref{lemma: dissipative A1} regarding the dissipativity of the operator $\mathcal{A}_1$. It will also be useful for assessing the exponential stability of the closed-loop system in the framework of energy multiplier method~\cite{Luo2012}.
\end{rem}
\subsection{$\mathcal{A}_1$ generates a $C_0$-semigroup of contractions}
We apply the Lumer-Phillips theorem~\cite{Curtain2012,Luo2012,Pazy2012} to show that $\mathcal{A}_1$ generates a $C_0$-semigroup. To do so, the following preliminary lemma is introduced.
\begin{lem} \label{lemma: inv A1}
The operator $\mathcal{A}_1^{-1} : \mathcal{H} \rightarrow D(\mathcal{A}_1)$ exists and is bounded, i.e., $\mathcal{A}_1^{-1}\in\mathcal{L}(\mathcal{H})$. Therefore, $0\in\rho(\mathcal{A}_1)$ and $\mathcal{A}_1$ is a closed operator.
\end{lem}
\textbf{Proof.}
Let us show first that $\mathcal{A}_1$ is surjective. Let $(\tilde{f},\tilde{g},\tilde{h},\tilde{z},\tilde{\zeta}_\omega,\tilde{\zeta}_\phi)\in\mathcal{H}$ be given. We are looking for $(f,g,h,z,\zeta_\omega,\zeta_\phi)\in D(\mathcal{A}_1)$ such that $\mathcal{A}_1(f,g,h,z,\zeta_\omega,\zeta_\phi)=(\tilde{f},\tilde{g},\tilde{h},\tilde{z},\tilde{\zeta}_\omega,\tilde{\zeta}_\phi)$. Introducing
\begin{equation*}
\begin{bmatrix}
\hat{g} \\ \hat{z}
\end{bmatrix}
\triangleq
\mathrm{diag}(1/\rho,1/I_w) M
\begin{bmatrix}
\tilde{g} \\ \tilde{z}
\end{bmatrix}
, \qquad
\begin{bmatrix} \hat{\zeta}_\omega \\ \hat{\zeta}_\phi \end{bmatrix}
\triangleq
M_s \begin{bmatrix} \tilde{\zeta}_\omega \\ \tilde{\zeta}_\phi \end{bmatrix} ,
\end{equation*}
it is equivalent to find $(f,g,h,z,\zeta_\omega,\zeta_\phi)\in D(\mathcal{A}_1)$ satisfying $g=\tilde{f}$, $z=\tilde{h}$, $\zeta_\omega = g(l) = \tilde{f}(l)$, $\zeta_\phi = z(l) = \tilde{h}(l)$,
\begin{subequations}
\begin{align}
-c_\omega^2 (f'' + \eta_\omega \tilde{f}'')'' = & \, \hat{g} , && \mathrm{in~}L^2(0,l) ; \label{eq: A1 surj 2}\\
c_\phi^2 (h'+\eta_\phi \tilde{h}')' = & \, \hat{z} , && \mathrm{in~}L^2(0,l) , \label{eq: A1 surj 4}
\end{align}
\end{subequations}
and the boundary conditions
\begin{subequations}
\begin{align}
\hat{\zeta}_\omega & = \rho c_\omega^2 (f''+\eta_\omega \tilde{f}'')'(l) - k_1 (\tilde{f}(l)+\epsilon_1 f(l)) , \label{eq: A1 surj - CL1} \\
\hat{\zeta}_\phi & = - I_w c_\phi^2 (h'+\eta_\phi \tilde{h}')(l) - k_2 (\tilde{h}(l)+\epsilon_2 h(l)) \label{eq: A1 surj - CL2} .
\end{align}
\end{subequations}
Direct computations show that
\begin{eqnarray}
f(y) & = & - \eta_\omega \tilde{f}(y) + \dfrac{k_1}{6} \cdot \dfrac{1 - \epsilon_1 \eta_\omega}{\rho c_\omega^2 + \epsilon_1 k_1 l^3/3} y^2(y-3l) \tilde{f}(l) \label{eq: surjectivite operator A1 def f} \\
& & + \dfrac{1}{6} \cdot \dfrac{1}{\rho c_\omega^2 + \epsilon_1 k_1 l^3/3} y^2(y-3l) \hat{\zeta}_\omega \nonumber \\
& &  - \dfrac{1}{c_\omega^2} \int_0^y\int_0^{\xi_1}\int_{\xi_2}^l\int_{\xi_3}^l \hat{g}(\xi_4) \mathrm{d}\xi_4\,\mathrm{d}\xi_3 \,\mathrm{d}\xi_2\,\mathrm{d}\xi_1 \nonumber \\
& &  - \dfrac{\epsilon_1 k_1}{6 c_\omega^2} \cdot \dfrac{1}{\rho c_\omega^2 + \epsilon_1 k_1 l^3/3} y^2(y-3l) \nonumber\\
& & \phantom{-\;} \times \int_0^l\int_0^{\xi_1}\int_{\xi_2}^l\int_{\xi_3}^l \hat{g}(\xi_4) \mathrm{d}\xi_4\,\mathrm{d}\xi_3\,\mathrm{d}\xi_2\,\mathrm{d}\xi_1 \nonumber  ,
\end{eqnarray}
and
\begin{eqnarray}
h(y) & = & - \eta_\phi \tilde{h}(y) - k_2 \dfrac{1 - \epsilon_2 \eta_\phi}{I_w c_\phi^2 + \epsilon_2 k_2 l} y \tilde{h}(l) \label{eq: surjectivite operator A1 def h} \\
& &  - \dfrac{1}{I_w c_\phi^2 + \epsilon_2 k_2 l} y \hat{\zeta}_\phi - \dfrac{1}{c_\phi^2} \int_0^y\int_{\xi_1}^l \hat{z}(\xi_2) \mathrm{d}\xi_2\,\mathrm{d}\xi_1 \nonumber \\
& &  + \dfrac{k_2}{c_\phi^2} \cdot \dfrac{\epsilon_2}{I_w c_\phi^2 + \epsilon_2 k_2 l} y \int_0^l\int_{\xi_1}^l \hat{z}(\xi_2) \mathrm{d}\xi_2\,\mathrm{d}\xi_1 \nonumber ,
\end{eqnarray}
solve (\ref{eq: A1 surj 2}-\ref{eq: A1 surj 4}) with $(f,g,h,z,\zeta_\phi,\zeta_\omega)\in D(\mathcal{A}_1)=D(\mathcal{A})$ while satisfying (\ref{eq: A1 surj - CL1}-\ref{eq: A1 surj - CL2}). Therefore,  $\mathcal{A}_1$ is surjective. Let us now investigate the injectivity. By definition, $\mathcal{A}_1(f,g,h,z,\zeta_\phi,\zeta_\omega)=(0,0,0,0,0,0)$ implies $g=z=0$ and thus $\zeta_\omega=g(l)=0$ and $\zeta_\phi=z(l)=0$. Hence, based on (\ref{eq: def A1 - part1}), $f''''=0$ and $h''=0$ in $L^2(0,l)$. As $(f,g,h,z,\zeta_\omega,\zeta_\phi)\in D(\mathcal{A}_1)$ and based on (\ref{eq: def A1 - part2}), the integration conditions are $f(0)=f'(0)=f''(l)=0$, $f'''(l)=(k_1 \epsilon_1 / (\rho c_\omega^2)) f(l)$, $h(0)=0$, and $h'(l) = -(k_2 \epsilon_2 / (I_w c_\phi^2)) h(l)$. Then, as $f''''(y)=0$ for almost all $y\in(0,l)$ and $f,f',f'',f'''\in\mathrm{AC}[0,l]$, it yields after four successive integrations that for any $y\in[0,l]$,
\begin{equation*}
f(y) = \dfrac{k_1 \epsilon_1}{\rho c_\omega^2} \cdot \dfrac{y^2(y-3l)}{6} f(l) .
\end{equation*}
Evaluating at $y=l$, it yields $(1+k_1 \epsilon_1 l^3 / (3 \rho c_\omega^2)) f(l) = 0$. As $l,\rho,c_\omega^2,k_1,\epsilon_1 \geq 0$, it implies that $f(l)=0$ and hence, $f=0$. Similarly, as $h''(y)=0$ for almost all $y\in(0,l)$ and $h,h'\in\mathrm{AC}[0,l]$, we have for any $y\in[0,l]$,
\begin{equation*}
h(y) = - \dfrac{k_2 \epsilon_2}{I_w c_\phi^2} y h(l) .
\end{equation*}
Evaluating at $y=l$, it yields $(1+k_2 \epsilon_2 l / (I_w c_\phi^2)) h(l) = 0$. As $l,I_w,c_\phi^2,k_2,\epsilon_2 \geq 0$, it implies that $h(l)=0$ and hence, $h=0$. Thus $\mathcal{A}_1(f,g,h,z,\zeta_\omega,\zeta_\phi)=(0,0,0,0,0,0)$ implies $(f,g,h,z,\zeta_\omega,\zeta_\phi)=(0,0,0,0,0,0)$ showing that $\mathcal{A}_1$ is injective.

Thus $\mathcal{A}_1 : D(\mathcal{A}_1) \rightarrow \mathcal{H}$ is bijective and $\mathcal{A}_1^{-1} : \mathcal{H} \rightarrow D(\mathcal{A}_1)$ is well defined for any $(\tilde{f},\tilde{g},\tilde{h},\tilde{z},\tilde{\zeta}_\omega,\tilde{\zeta}_\phi)\in\mathcal{H}$ by $\mathcal{A}_1^{-1}(\tilde{f},\tilde{g},\tilde{h},\tilde{z},\tilde{\zeta}_\omega,\tilde{\zeta}_\phi) = (f,\tilde{f},h,\tilde{h},\tilde{f}(l),\tilde{h}(l))$ where $f$ and $h$ are given by equations (\ref{eq: surjectivite operator A1 def f}) and (\ref{eq: surjectivite operator A1 def h}), respectively. Finally, by employing Schwartz's and Poincar{\'e}'s inequalities, it is straightforward to show that $\mathcal{A}_1^{-1}$ is bounded, i.e., there exists $C \in\mathbb{R}_+$ such that for all $X\in\mathcal{H}$, $\norme{\mathcal{A}_1^{-1}X}_{\mathcal{H},1} \leq C \norme{X}_{\mathcal{H},1}$. It shows that $\mathcal{A}_1^{-1}\in\mathcal{L}(\mathcal{H})$. Then $0\in\rho(\mathcal{A}_1)$ and $\mathcal{A}_1$ is a closed operator.
\hfill $\qed$

The second key-element for applying the Lumer-Phillips theorem is stated in the following lemma with $\epsilon_1^*,\epsilon_2^* \in \mathbb{R}_+^*$ defined by
\begin{align*}
\epsilon_1^* = & \min \left(
\dfrac{\pi^2 I_w \eta_\phi c_\phi^2}{l |x_c| (4 l \rho + \pi^2 m_s )} , \right. \\
& \left. \dfrac{2 \pi^4 \rho \eta_\omega c_\omega^2}{32 l^4 (2+|x_c|) \rho + 8 \pi^2 l^3 (2+|x_c|) m_s + \pi^4 \rho \eta_\omega^2 c_\omega^2}
\right) , \\
\epsilon_2^* = & \min \left(
\dfrac{\pi^4 \rho \eta_\omega c_\omega^2}{4 l^3 |x_c| (4 l \rho + \pi^2 m_s)} , \right. \\
& \left. \dfrac{2 \pi^2 I_w \eta_\phi c_\phi^2}{8 l^2 (2 I_w^* + |x_c| \rho) + 2 \pi^2 l (2 J_s^* + |x_c| m_s) + \pi^2 I_w \eta_\phi^2 c_\phi^2}
\right) .
\end{align*}
\begin{lem} \label{lemma: dissipative A1}
Let $\epsilon_1,\epsilon_2 \in \mathbb{R}_+^*$ such that $\epsilon_1 < \min( \epsilon_1^* , K_{m,1})$ and $\epsilon_2 < \min( \epsilon_2^* , K_{m_2})$. Then the operator $\mathcal{A}_1 : D(\mathcal{A}_1) \rightarrow \mathcal{H}$ is dissipative with respect to $\left\langle  \cdot , \cdot \right\rangle_{\mathcal{H},2}$.
\end{lem}
\textbf{Proof.} Lemma~\ref{lemma: inner product 2} ensures that $(\mathcal{H},\left\langle  \cdot , \cdot \right\rangle_{\mathcal{H},2})$ is a real Hilbert space. Thus, to prove the dissipativity of operator $\mathcal{A}_1$, we have to show that $\left\langle \mathcal{A}_1 X , X \right\rangle_{\mathcal{H},2} \leq 0$ for all $X \in D(\mathcal{A}_1)$~\cite{Curtain2012,Luo2012,Pazy2012}. Let $X=(f,g,h,z,\zeta_\omega,\zeta_\phi) \in D(\mathcal{A}_1)$. Straightforward integrations by parts for absolutely continuous functions~\cite{Cannarsa2015} along with the boundary conditions (\ref{eq: def domain A}) yields :
\begin{align}
\langle \mathcal{A}_1 & X , X \rangle_{\mathcal{H},2} \nonumber \\
= & - k_1 (g(l)+\epsilon_1 f(l))^2 - k_2 (z(l)+\epsilon_2 h(l))^2 \label{eq: inner A1X-X IPP}  \\
& + \epsilon_1 \rho \int_0^l g(y)^2 \mathrm{d}y + \epsilon_2 I_w^* \int_0^l z(y)^2 \mathrm{d}y \nonumber \\
& + (\epsilon_1 + \epsilon_2) \rho x_c \int_0^l g z \mathrm{d}y \nonumber \\
& - \rho \eta_\omega c_\omega^2 \int_0^l g''(y)^2 \mathrm{d}y - I_w \eta_\phi c_\phi^2 \int_0^l z'(y)^2 \mathrm{d}y \nonumber \\
& - \epsilon_1 \rho c_\omega^2 \int_0^l f''(y)^2 \mathrm{d}y - \epsilon_1 \rho \eta_\omega c_\omega^2 \int_0^l f''(y)  g''(y) \mathrm{d}y \nonumber \\
& - \epsilon_2 I_w c_\phi^2 \int_0^l h'(y)^2 \mathrm{d}y - \epsilon_2 I_w \eta_\phi c_\phi^2 \int_0^l h'(y)  z'(y) \mathrm{d}y \nonumber \\
& + \epsilon_1 m_s g(l)^2 + \epsilon_2 J_s^* z(l)^2 + (\epsilon_1+\epsilon_2) m_s x_c g(l) z(l). \nonumber
\end{align}
Applying Young's inequality with $r_1,r_2>0$ to be determined later, it provides
\begin{eqnarray*}
\left| \int_0^l f''(y) g''(y) \mathrm{d}y \right| & \leq & \dfrac{1}{2 r_1} \int_0^l f''(y)^2 \mathrm{d}y + \dfrac{r_1}{2} \int_0^l g''(y)^2 \mathrm{d}y , \\
\left| \int_0^l h'(y) z'(y) \mathrm{d}y \right| & \leq & \dfrac{1}{2 r_2} \int_0^l h'(y)^2 \mathrm{d}y + \dfrac{r_2}{2} \int_0^l z'(y)^2 \mathrm{d}y .
\end{eqnarray*}
Furthermore, the terms in the last line of (\ref{eq: inner A1X-X IPP}) are handled by using Schwartz's inequality: $g(l)^2 = \left( \int_0^l g'(y) \mathrm{dy} \right)^2 \leq l \int_0^l g'(y)^2 \mathrm{d}y$ and, similarly, $z(l)^2 \leq l \int_0^l z'(y)^2 \mathrm{d}y$. Finally, by resorting to Poincar\'{e}'s inequality, it yields for any $X=(f,g,h,z,\zeta_\omega,\zeta_\phi)\in D(\mathcal{A}_1)$ and any $r_1,r_2>0$,
\begin{eqnarray}
& & \left\langle \mathcal{A}_1 X , X \right\rangle_{\mathcal{H},2} \nonumber \\
& \leq & - k_1 (g(l)+\epsilon_1 f(l))^2 - k_2 (z(l)+\epsilon_2 h(l))^2 \label{eq: inner A1X-X} \\
& & - \rho \nu_1 \int_0^l g''(y)^2 \mathrm{d}y \nonumber - I_w \nu_2 \int_0^l z'(y)^2 \mathrm{d}y \nonumber \\
& & - \epsilon_1 \left( 1 - \dfrac{\eta_\omega}{2 r_1} \right) \int_0^l \rho c_\omega^2 f''(y)^2 \mathrm{d}y \nonumber \\
& & - \epsilon_2 \left( 1 - \dfrac{\eta_\phi}{2 r_2} \right) \int_0^l I_w c_\phi^2 h'(y)^2 \mathrm{d}y . \nonumber
\end{eqnarray}
where
\begin{subequations}
\begin{align}
\nu_1 \triangleq & \eta_\omega c_\omega^2 - \epsilon_2 \dfrac{2 l^3 |x_c|}{\pi^2} \left( \dfrac{4 l}{\pi^2} + \dfrac{m_s}{\rho} \right) \label{eq: nu1} \\
& - \epsilon_1 \left( \dfrac{8 l^4 (2+|x_c|)}{\pi^4} + \dfrac{\eta_\omega c_\omega^2 r_1}{2} + \dfrac{2 l^3 (2+|x_c|) m_s}{\pi^2\rho} \right) , \nonumber \\
\nu_2 \triangleq & \eta_\phi c_\phi^2 - \epsilon_1 \dfrac{l |x_c|}{I_w} \left( \dfrac{2 l \rho}{\pi^2} + \dfrac{m_s}{2} \right) \label{eq: nu2} \\
& - \epsilon_2 \left( 2 l^2 \dfrac{2 I_w^* + |x_c| \rho}{\pi^2 I_w} + \dfrac{\eta_\phi c_\phi^2 r_2}{2} + \dfrac{l (2 J_s^* + |x_c| m_s)}{2 I_w} \right) . \nonumber
\end{align}
\end{subequations}
As $\epsilon_2 < \epsilon_2^*$,
\begin{align*}
\nu_1 \geq & \dfrac{\eta_\omega c_\omega^2}{2} - \epsilon_1 \left( \dfrac{8 l^4 (2+|x_c|)}{\pi^4} + \dfrac{\eta_\omega c_\omega^2 r_1}{2} + \dfrac{2 l^3 (2+|x_c|) m_s}{\pi^2\rho} \right) \\
\geq &\left( \dfrac{8 l^4 (2+|x_c|)}{\pi^4} + \dfrac{\eta_\omega c_\omega^2 r_1}{2} + \dfrac{2 l^3 (2+|x_c|) m_s}{\pi^2\rho} \right) \\
     & \times(\varphi_1(r_1) - \epsilon_1) ,
\end{align*}
where 
\begin{equation*}
\varphi_1(x) = \dfrac{\pi^4\rho\eta_\omega c_\omega^2}{16 l^4 (2+|x_c|) \rho + 4 \pi^2 l^3 (2+|x_c|) m_s + \pi^4 \rho \eta_\omega c_\omega^2 x}
\end{equation*}
is a continuous decreasing function over $\mathbb{R}_+^*$ that tends to zero when $x \rightarrow +\infty$. By assumption we have $\epsilon_1 < \epsilon_1^* \leq \varphi_1(\eta_\omega/2)$. Hence there exists $r_1^* > \eta_\omega /2$ such that $\epsilon_1 < \varphi_1(r_1^*) < \varphi_1(\eta_\omega/2)$, implying $\nu_1 > 0$ for $r_1=r_1^*$. Similarly one can show that there exists $r_2^*>\eta_\phi/2$ such that $\nu_2 > 0$ for $r_2=r_2^*$. Therefore, taking $r_1=r_1^*$ and $r_2=r_2^*$ in (\ref{eq: inner A1X-X}), it ensures that for any $X\in\mathcal{H}$, $\left\langle \mathcal{A}_1 X , X \right\rangle_{\mathcal{H},2} \leq 0$, i.e., $\mathcal{A}_1$ is dissipative on $\mathcal{H}$ endowed with $\left\langle  \cdot , \cdot \right\rangle_{\mathcal{H},2}$.
\hfill $\qed$

We are now ready to establish the main property of the operator $\mathcal{A}_1$.
\begin{thm}\label{th: A1 generation semigroup}
Let $\epsilon_1,\epsilon_2 \in \mathbb{R}_+^*$ such that $\epsilon_1 < \min( \epsilon_1^* , K_{m,1})$ and $\epsilon_2 < \min( \epsilon_2^* , K_{m,2})$. Then the operator $\mathcal{A}_1$ generates a $C_0$-semigroup of contractions on $(\mathcal{H},\left\langle  \cdot , \cdot \right\rangle_{\mathcal{H},2})$.
\end{thm}
\textbf{Proof.} Lemma~\ref{lemma: inner product 2} ensures that $(\mathcal{H},\left\langle  \cdot , \cdot \right\rangle_{\mathcal{H},2})$ is a real Hilbert space. Furthermore, Lemma~\ref{lemma: dissipative A1} shows that the linear operator $\mathcal{A}_1$ is dissipative with respect to $\left\langle  \cdot , \cdot \right\rangle_{\mathcal{H},2}$. It remains to show that there exists $\lambda_0 > 0 $ such that $R(\lambda_0 I_{D(\mathcal{A}_1)} - \mathcal{A}_1) = \mathcal{H}$. Lemma~\ref{lemma: inv A1} shows that the operator $\mathcal{A}_1$ is closed. Then its resolvent set $\rho(\mathcal{A}_1)$ is an open subset of $\mathbb{C}$. As $0\in\rho(\mathcal{A}_1)$, there exists $\lambda_0>0$ such that the range condition $R(\lambda_0 I_{D(\mathcal{A}_1)} - \mathcal{A}_1) = \mathcal{H}$ holds. The application of the Lumer-Philips theorem~\cite[Th.2.29]{Luo2012}~\cite{Pazy2012} for reflexive spaces concludes the proof.
\hfill $\qed$
\begin{cor}\label{cor: density}
Let $\epsilon_1,\epsilon_2 \in \mathbb{R}_+^*$ such that $\epsilon_1 < \min( \epsilon_1^* , K_{m,1})$ and $\epsilon_2 < \min( \epsilon_2^* , K_{m,2})$. Then, $D(\mathcal{A})$ is dense in $\mathcal{H}$ endowed by either $\langle\cdot,\cdot\rangle_{\mathcal{H},1}$ or $\langle\cdot,\cdot\rangle_{\mathcal{H},2}$, i.e., $\overline{D(\mathcal{A})}=\mathcal{H}$.
\end{cor}
\textbf{Proof.} The property that $D(\mathcal{A}_1)$ is dense in $\mathcal{H}$ endowed by $\langle\cdot,\cdot\rangle_{\mathcal{H},2}$ follows from both dissipativity and range conditions satisfied by operator $\mathcal{A}_1$ in the application of the Lumer-Philips theorem for reflexive spaces~\cite[Chap.1, Th.4.5 and Th.4.6]{Pazy2012}. As the norms $\norme{\cdot}_{\mathcal{H},1}$ and $\norme{\cdot}_{\mathcal{H},2}$ are equivalent, $D(\mathcal{A}_1)$ is also dense in $\mathcal{H}$ endowed by $\langle\cdot,\cdot\rangle_{\mathcal{H},1}$. The proof is complete because, by definition, $D(\mathcal{A})=D(\mathcal{A}_1)$. \hfill $\qed$

Based on (\ref{eq: inner A1X-X}) and under the assumptions of Theorem~\ref{th: A1 generation semigroup}, it can be shown that the $C_0$-semigroup generated by $\mathcal{A}_1$ is exponentially stable (the proof follows the one of Theorem~\ref{th: exponential stab C0-semigroup 2nd inner product}). Thus, assuming as in~\cite{He2017,Paranjape2013} bounded aerodynamic efforts $F_a$ and $M_a$, it can be concluded that both the system energy and the flexible displacements are bounded. In this work, adopting the aerodynamic model (\ref{eq: aero efforts}), we derive a sufficient condition for ensuring the exponential stability of (\ref{eq: cauchy problem}).

\subsection{$\mathcal{A}$ generates a $C_0$-semigroup}
We can now introduce the main result of this section.
\begin{thm}
Let $\epsilon_1,\epsilon_2 \in \mathbb{R}_+^*$ such that $\epsilon_1 < \min( \epsilon_1^* , K_{m,1})$ and $\epsilon_2 < \min( \epsilon_2^* , K_{m,2})$. Then the operator $\mathcal{A}$ generates a $C_0$-semigroup on both $(\mathcal{H},\left\langle  \cdot , \cdot \right\rangle_{\mathcal{H},1})$ and $(\mathcal{H},\left\langle  \cdot , \cdot \right\rangle_{\mathcal{H},2})$.
\end{thm}
\textbf{Proof.} We have by definition $\mathcal{A}=\mathcal{A}_1+\mathcal{A}_2$ over $D(\mathcal{A})$ with $D(\mathcal{A})=D(\mathcal{A}_1)$ and $D(\mathcal{A}_2)=\mathcal{H}$. By Theorem~\ref{th: A1 generation semigroup}, $\mathcal{A}_1$ generates a $C_0$-semigroup of contractions with respect to $\left\langle  \cdot , \cdot \right\rangle_{\mathcal{H},2}$. Moreover, it can be shown from its definition that $\mathcal{A}_2$ is a bounded operator with respect to $\left\Vert  \cdot \right\Vert_{\mathcal{H},2}$. Then the perturbation theory~\cite[Chap.3, Th.1.1]{Pazy2012}\cite[Th.3.2.1]{Curtain2012} ensures that $\mathcal{A}$ generates a $\mathcal{C}_0$-semigroup on $(\mathcal{H},\left\langle  \cdot , \cdot \right\rangle_{\mathcal{H},2})$ and hence on $(\mathcal{H},\left\langle  \cdot , \cdot \right\rangle_{\mathcal{H},1})$ based on Lemma~\ref{lemma: inner product 2}.
\hfill $\qed$

Therefore, the Cauchy problem (\ref{eq: cauchy problem}) is well-posed~\cite{Curtain2012,Luo2012,Pazy2012}. Let $T:\mathbb{R}_+ \rightarrow \mathcal{L}(\mathcal{H})$ be the $C_0$-semigroup generated by $\mathcal{A}$. Then, for any initial condition $X_0 \in D(\mathcal{A})$, $X(t)\triangleq T(t)X_0$ defined for $t \geq 0$ is the unique solution of (\ref{eq: cauchy problem}). Furthermore, for any $t\geq 0$, $X(t)\in D(\mathcal{A})$~\cite[Chap.1, Th.2.4]{Pazy2012}, i.e., $X\in\mathcal{C}^0(\mathbb{R}_+;D(\mathcal{A})) \cap \mathcal{C}^1(\mathbb{R}_+;\mathcal{H})$.

\section{Exponential Stability Assessment}\label{sec: stability}

\subsection{Exponential decay of the system energy}
We introduce, in the framework of energy multiplier methods~\cite{Luo2012}, the augmented energy of the plant, defined by
\begin{equation}\label{eq: definition augmented energy}
\forall t \geq 0, \; \mathcal{E}(t) \triangleq \dfrac{1}{2} \norme{X(t)}_{\mathcal{H},2}^2 = \dfrac{1}{2} \left\langle X(t) , X(t) \right\rangle_{\mathcal{H},2} ,
\end{equation}
where $X(t)\triangleq T(t) X_0 \in D(\mathcal{A})$ denotes the unique solution of (\ref{eq: cauchy problem}) associated to the initial condition $X_0\in D(\mathcal{A})$. It can be shown that $\mathcal{E}$ is continuously differentiable over $\mathbb{R}_+$, with a derivative $\dot{\mathcal{E}}$ satisfying for any $t \geq 0$:
\begin{eqnarray}
\dot{\mathcal{E}}(t) & = & \left\langle \dot{X}(t) , X(t) \right\rangle_{\mathcal{H},2} = \left\langle \mathcal{A}X(t) , X(t) \right\rangle_{\mathcal{H},2} \nonumber \\
& = & \left\langle \mathcal{A}_1X(t) , X(t) \right\rangle_{\mathcal{H},2} + \left\langle \mathcal{A}_2X(t) , X(t) \right\rangle_{\mathcal{H},2} . \label{eq: derivative energie E}
\end{eqnarray}

In order to derive the exponential stability of the system, we formulate the following assumption.
\begin{assum}\label{assumption: physical parameters}
It is assumed that the constants in (\ref{eq: studied EDP}) are such that there exist $r_1,r_2,\ldots,r_8>0$ along with $0 < \epsilon_1 < \min( \epsilon_1^* , K_{m,1})$ and $0 < \epsilon_2 < \min( \epsilon_2^* , K_{m,2})$ such that
\begin{eqnarray*}
\lambda_1 & \triangleq & \epsilon_1 \left( 1 - \dfrac{\eta_\omega}{2 r_1} - \dfrac{8 l^4}{\pi^4 \rho c_\omega^2} \left( \dfrac{\alpha_\omega}{r_4} + \dfrac{\beta_\omega}{r_5} + \dfrac{\gamma_\omega}{r_3} \right) \right) , \\
\lambda_2 & \triangleq & \gamma_\omega + \dfrac{\alpha_\omega + \epsilon_2 \gamma_\phi}{2 r_6} + \dfrac{\epsilon_1 \gamma_\omega r_3}{2} + \dfrac{\beta_\omega + \gamma_\phi}{2 r_7} \geq 0, \\
\lambda_3 & \triangleq & \rho \left\{ \eta_\omega c_\omega^2 - \epsilon_2 \dfrac{2 l^3 |x_c|}{\pi^2} \left( \dfrac{4 l}{\pi^2} + \dfrac{m_s}{\rho} \right) \right. \\
& & \phantom{=} \left. - \epsilon_1 \left( \dfrac{8 l^4 (2+|x_c|)}{\pi^4} + \dfrac{\eta_\omega c_\omega^2 r_1}{2} + \dfrac{2 l^3 (2+|x_c|) m_s}{\pi^2\rho} \right) \right\} , \\
\lambda_4 & \triangleq & \epsilon_2 \left( 1 - \dfrac{\eta_\phi}{2 r_2} \right) - \dfrac{4 l^2}{\pi^2 I_w c_\phi^2} \left( \dfrac{(\alpha_\omega + \epsilon_2 \gamma_\phi)r_6}{2} + \dfrac{\alpha_\phi + \epsilon_2 \beta_\phi}{2 r_8} \right. \\
& & \phantom{=================} \left. + \dfrac{\epsilon_1 \alpha_\omega r_4}{2} + \epsilon_2 \alpha_\phi \right) , \\
\lambda_5 & \triangleq & \beta_\phi + \dfrac{(\beta_\omega + \gamma_\phi)r_7}{2} + \dfrac{(\alpha_\phi + \epsilon_2 \beta_\phi)r_8}{2} + \dfrac{\epsilon_1 \beta_\omega r_5}{2} \geq 0 , \\
\lambda_6 & \triangleq & I_w \left\{ \eta_\phi c_\phi^2 - \epsilon_1 \dfrac{l |x_c|}{I_w} \left( \dfrac{2 l \rho}{\pi^2} + \dfrac{m_s}{2} \right) \right. \\
& & \phantom{==} \left. - \epsilon_2 \left( 2 l^2 \dfrac{2 I_w^* + |x_c| \rho}{\pi^2 I_w} + \dfrac{\eta_\phi c_\phi^2 r_2}{2} + \dfrac{l (2 J_s^* + |x_c| m_s)}{2 I_w} \right) \right\} ,
\end{eqnarray*}
satisfy $\lambda_1, \lambda_3, \lambda_4 , \lambda_6 > 0$, $\mu_3 \triangleq \pi^4 \lambda_3 / (16 l^4) - \lambda_2 > 0$, and $\mu_6 \triangleq \pi^2 \lambda_6 / (4l^2) - \lambda_5 > 0$.
\end{assum}
\begin{rem}
Similarly to~\cite{Bialy2016,He2017}, Assumption~\ref{assumption: physical parameters} imposes restrictions on both control parameters and the physical parameters of the wing. It is essentially a trade-off between the structural stiffness and the amplitude of the aerodynamic coefficients. In particular, it can be seen that the constraints can always be met by increasing the stiffness parameters. Indeed, for fixed aerodynamic coefficients,
\begin{equation*}
\lambda_1 \underset{c_\omega \rightarrow \infty}{\sim} \epsilon_1 \left( 1 - \dfrac{\eta_\omega}{2 r_1} \right)
, \qquad
\lambda_4 \underset{c_\phi \rightarrow \infty}{\sim} \epsilon_2 \left( 1 - \dfrac{\eta_\phi}{2 r_2} \right) ,
\end{equation*}
\begin{equation*}
\mu_3 \underset{c_\omega \rightarrow \infty}{\sim} \pi^4 \lambda_3 / (16 l^4)
, \qquad
\mu_6 \underset{c_\phi \rightarrow \infty}{\sim} \pi^2 \lambda_6 / (4l^2) .
\end{equation*}
Noting that $\lambda_3 = \rho \nu_1$ and $\lambda_6 = I_w \nu_2$ (see (\ref{eq: nu1})-(\ref{eq: nu2})), it corresponds to the constraints involved in the proof of Lemma~\ref{lemma: dissipative A1} for which it was shown that a feasible solution exists for any arbitrary value of the structural parameters. Thus, while the exponential stability of the closed-loop aerodynamic free model is always ensured, the one of the full model requires an adequate balance between the structural stiffness and the aerodynamic parameters. To fulfill the constraints, one can resort to the following design procedure. First, aerodynamic coefficients are obtained based on performance criteria corresponding to the desired flight envelope. Then, structural stiffness, along with the control parameters $\epsilon_1,\epsilon_2$ can be adjusted to satisfy Assumption~\ref{assumption: physical parameters}.
\end{rem}

\begin{thm}\label{th: exponential stab C0-semigroup 2nd inner product}
Assume that Assumption~\ref{assumption: physical parameters} holds. Then the augmented energy $\mathcal{E}$ defined by (\ref{eq: definition augmented energy}) exponentially decays to zero, i.e., there exists $\Lambda \in\mathbb{R}_+^*$ such that
\begin{equation*}
\forall t \geq 0,  \;
\mathcal{E}(t) \leq \mathcal{E}(0) \exp(-\Lambda t).
\end{equation*}
Furthermore, $T(t)$ is an exponentially stable $C_0$-semigroup for $\norme{\cdot}_{\mathcal{H},2}$.
\end{thm}
\textbf{Proof.}
As $\mathcal{E}\in\mathcal{C}^1(\mathbb{R}_+;\mathbb{R})$, the objective is to show that there exists a $\Lambda\in\mathbb{R}_+^*$ such that for any $t\geq0$, $\dot{\mathcal{E}}(t) \leq - \Lambda \mathcal{E}(t)$. We first note that the first term of the right-hand side of (\ref{eq: derivative energie E}) is upper bounded by (\ref{eq: inner A1X-X}) since for any $t\geq 0$, $X(t)\in D(\mathcal{A})$ and, by definition of operator $\mathcal{A}_1$, $D(\mathcal{A}_1)=D(\mathcal{A})$. We study the second term in the right-hand side of (\ref{eq: derivative energie E}). For any $X = (f,g,h,z,\zeta_\omega,\zeta_\phi)\in D(\mathcal{A})\subset D(\mathcal{A}_2)$, applying first Young's inequality and then Poincar\'{e}'s inequality, it yields for any $r_3,\ldots,r_8>0$, and in particular for the value of the parameters satisfying Assumption~\ref{assumption: physical parameters},
\begin{eqnarray*}
& & \left\langle \mathcal{A}_2X , X \right\rangle_{\mathcal{H},2} \\
& = & \int_0^l F_a(y) g(y) + M_a(y) z(y) + \epsilon_1 F_a(y) f(y) + \epsilon_2 M_a(y) h(y) \mathrm{d}y \\
& \leq & \dfrac{8 l^4 \epsilon_1}{\pi^4 \rho c_\omega^2} \left( \dfrac{\alpha_\omega}{r_4} + \dfrac{\beta_\omega}{r_5} + \dfrac{\gamma_\omega}{r_3} \right) \int_0^{l} \rho c_\omega^2 f''(y)^2 \mathrm{d}y \\
&  & +  \left( \gamma_\omega + \dfrac{\alpha_\omega + \epsilon_2 \gamma_\phi}{2 r_6} + \dfrac{\epsilon_1 \gamma_\omega r_3}{2} + \dfrac{\beta_\omega + \gamma_\phi}{2 r_7} \right) \int_0^{l} g(y)^2 \mathrm{d}y \\
&  & +  \dfrac{4 l^2}{\pi^2 I_w c_\phi^2} \left( \dfrac{(\alpha_\omega + \epsilon_2 \gamma_\phi)r_6}{2} + \dfrac{\alpha_\phi + \epsilon_2 \beta_\phi}{2 r_8} + \dfrac{\epsilon_1 \alpha_\omega r_4}{2} + \epsilon_2 \alpha_\phi \right) \\
& & \phantom{+\;} \times \int_0^{l} I_w c_\phi^2 h'(y)^2 \mathrm{d}y \\
&  & +  \left( \beta_\phi + \dfrac{(\beta_\omega + \gamma_\phi)r_7}{2} + \dfrac{(\alpha_\phi + \epsilon_2 \beta_\phi)r_8}{2} + \dfrac{\epsilon_1 \beta_\omega r_5}{2} \right) \\
& & \phantom{+\;} \times \int_0^{l} z(y)^2 \mathrm{d}y .
\end{eqnarray*}
Let $X(t) = (f(\cdot,t),g(\cdot,t),h(\cdot,t),z(\cdot,t),\zeta_\omega(t),\zeta_\phi(t))\in D(\mathcal{A})$ be the solution to (\ref{eq: cauchy problem}) associated to the initial condition $X(0)=X_0\in D(\mathcal{A})$. Considering the values of $r_1,\ldots,r_8,\epsilon_1,\epsilon_2$ satisfying Assumption~\ref{assumption: physical parameters}, together with (\ref{eq: inner A1X-X}), (\ref{eq: derivative energie E}) verifies for any $t\geq0$,
\begin{eqnarray*}
\dot{\mathcal{E}}(t)
& \leq & - \lambda_1 \int_0^l \rho c_\omega^2 f''(y,t)^2 \mathrm{d}y + \lambda_2  \int_0^{l} g(y,t)^2 \mathrm{d}y \\
& & - \lambda_3 \int_0^l g''(y,t)^2 \mathrm{d}y - \lambda_4 \int_0^l I_w c_\phi^2 h'(y,t)^2 \mathrm{d}y \\
& & + \lambda_5 \int_0^{l} z(y,t)^2 \mathrm{d}y - \lambda_6 \int_0^l z'(y,t)^2 \mathrm{d}y .
\end{eqnarray*}
As $\mu_3>0$, there exists $\lambda_3^*\in(0,\lambda_3)$ such that $\Delta\lambda_3 \triangleq \lambda_3 - \lambda_3^* > 0$ and $\mu_3^* \triangleq \pi^4 \lambda_3^* / (16 l^4) - \lambda_2 > 0$. Similarly, there exists $\lambda_6^*\in(0,\lambda_6)$ such that $\Delta\lambda_6 \triangleq \lambda_6 - \lambda_6^* > 0$ and $\mu_6^* \triangleq \pi^2 \lambda_6^* / (4 l^2) - \lambda_5 > 0$. Then, based on Poincar\'{e}'s inequality,
\begin{eqnarray}
\dot{\mathcal{E}}(t)
& \leq & - \lambda_1 \int_0^l \rho c_\omega^2 f''(y,t)^2 \mathrm{d}y - \mu_3^*  \int_0^{l} g(y,t)^2 \mathrm{d}y \label{eq: upper bound dot_mathcal_E} \\
& & - \lambda_4 \int_0^l I_w c_\phi^2 h'(y,t)^2 \mathrm{d}y - \mu_6^* \int_0^{l} z(y,t)^2 \mathrm{d}y \nonumber \\
& & - \Delta\lambda_3 \int_0^l g''(y,t)^2 \mathrm{d}y - \Delta\lambda_6 \int_0^l z'(y,t)^2 \mathrm{d}y \nonumber .
\end{eqnarray}
By applying Schwartz's inequality, we have $\zeta_\phi(t)^2 = z(l,t)^2 = \left(\int_0^l z'(y,t) \mathrm{d}y\right)^2 \leq l \Vert z'(\cdot,t) \Vert_{L^2(0,l)}^2$. Similarly, one can get $\zeta_\omega(t)^2 \leq (4l^3/\pi^2) \Vert g''(\cdot,t) \Vert_{L^2(0,l)}^2$. We then introduce $\mu_m>0$ defined by
\begin{equation*}
\mu_m \triangleq 2 \min\left( \lambda_1 , \lambda_4, \dfrac{\mu_3^*}{\lambda_M(M)}, \dfrac{\mu_6^*}{\lambda_M(M)} , \dfrac{\pi^2 \Delta\lambda_3}{4 l^3 \lambda_M(M_s)} , \dfrac{\Delta\lambda_6}{l \lambda_M(M_s)} \right) ,
\end{equation*}
where $\lambda_M(\cdot)$ denotes the largest eigenvalue. Recalling that for any symmetric matrix $S\in\mathbb{R}^{n \times n}$ and any vector $x\in\mathbb{R}^n$, $x^\top S x \leq \lambda_M(S) x^\top x$, it provides, based first on (\ref{eq: upper bound dot_mathcal_E}) and then on (\ref{eq: norm equivalent}),
\begin{equation*}
\forall t\geq0 , \; \dot{\mathcal{E}}(t) \leq - \dfrac{\mu_m}{2} \norme{X(t)}_{\mathcal{H},1}^2 \leq - \dfrac{\Lambda}{2} \norme{X(t)}_{\mathcal{H},2}^2 = - \Lambda \mathcal{E}(t) ,
\end{equation*}
where $\Lambda \triangleq \mu_m/(1+\alpha)>0$ is independent of the initial condition $X_0\in D(\mathcal{A})$. Thus, for any $t \geq 0$, $\mathcal{E}(t) \leq \mathcal{E}(0) \exp(-\Lambda t)$. Recalling that $\mathcal{E}(t)=\norme{X(t)}_{\mathcal{H},2}^2/2$, it yields that
\begin{equation*}
\forall X_0 \in D(\mathcal{A}),\; \forall t\geq0 ,\; \norme{T(t)X_0}_{\mathcal{H},2} \leq \norme{X_0}_{\mathcal{H},2} \exp(-\Lambda t/2) .
\end{equation*}
This inequality can be extended by density to all $X_0\in\mathcal{H}$ because $T(t)\in\mathcal{L}(\mathcal{H})$ and, based on Corollary~\ref{cor: density}, $\overline{D(\mathcal{A})}=\mathcal{H}$. Thus $T(t)$ is an exponentially stable $C_0$-semigroup for $\norme{\cdot}_{\mathcal{H},2}$ with $\left|\left|\left|T(t)\right|\right|\right|_{\mathcal{H},2} \leq \exp(-\Lambda t/2)$. \hfill $\qed$
\begin{cor}
Assume that Assumption~\ref{assumption: physical parameters} holds. Then the energy $E$ defined by (\ref{eq: definition energy}) exponentially decays to zero. In particular, there exists $K_E\in\mathbb{R}_+^*$ such that
\begin{equation*}
\forall t \geq 0,  \;
E(t) \leq K_E E(0) \exp(-\Lambda t) ,
\end{equation*}
where $\Lambda>0$ is provided by Theorem~\ref{th: exponential stab C0-semigroup 2nd inner product}. Furthermore, $T(t)$ is an exponentially stable $C_0$-semigroup for $\norme{\cdot}_{\mathcal{H},1}$.
\end{cor}
\textbf{Proof.} With $K_E = (1 + \alpha)/(1 - \alpha)>0$, it is a direct consequence of the equivalence of the norms (\ref{eq: norm equivalent}).
\hfill $\qed$

\subsection{Uniform exponential stability of the bending and twisting displacements}
Finally, we assess the uniform exponential stability of bending and twisting displacements.
\begin{cor}\label{cor: uniform exp decay flexible displacements}
Assume that Assumption~\ref{assumption: physical parameters} holds. Then, there exist $K_\omega,K_{\omega_y},K_\phi\in\mathbb{R}_{+}^{*}$ such that for any initial condition $X_0\in D(\mathcal{A})$ and for all $t\geq0$, the solution to the Cauchy problem (\ref{eq: cauchy problem}) satisfies for $\psi \in \{\omega,\omega_y,\phi\}$,
\begin{equation}\label{eq: uniform exp decay to zero of flexible displacements}
\norme{\psi(\cdot,t)}_\infty \leq K_\psi \sqrt{E(0)} \exp(-\Lambda t /2) ,
\end{equation}
where $\Lambda>0$ is given in Theorem~\ref{th: exponential stab C0-semigroup 2nd inner product} and $E(0)=\norme{X_0}_{\mathcal{H},1}^2/2$ is the initial energy of the system.
\end{cor}
\textbf{Proof.} Let $X(t) = (f(\cdot,t),g(\cdot,t),h(\cdot,t),z(\cdot,t),\zeta_\omega(t),\zeta_\phi(t))$ be the solution of the Cauchy problem (\ref{eq: cauchy problem}) associated to the initial condition $X_0\in D(\mathcal{A})$. As $X(t)\in D(\mathcal{A})$ for all $t \geq 0$, Agmon's inequality and then Poincar\'{e}'s inequality yield
\begin{eqnarray*}
\left|\left|h(\cdot,t)\right|\right|_\infty ^ 4
& \leq & 4 \int_0^l h(y,t)^2 \mathrm{d}y \int_0^l h'(y,t)^2 \mathrm{d}y \\
& \leq & \dfrac{16 l^2}{\pi^2} \left[ \int_0^l h'(y,t)^2 \mathrm{d}y \right]^2 \leq \dfrac{64 l^2}{\pi^2 I_w^2 c_\phi^4} E(t)^2 .\\
\end{eqnarray*}
Hence, as $h=\phi$, (\ref{eq: uniform exp decay to zero of flexible displacements}) holds with
\begin{equation*}
K_\phi = \dfrac{2}{c_\phi} \sqrt{\dfrac{2 l}{\pi I_w} \cdot \dfrac{1 + \alpha}{1 - \alpha}} .
\end{equation*}
Following the same procedure, similar inequalities are obtained for $f=\omega$ and $f'=\omega_y$. \hfill $\qed$

We deduce from Theorem~\ref{th: exponential stab C0-semigroup 2nd inner product} and Corollary~\ref{cor: uniform exp decay flexible displacements} that the command given by (\ref{eq: bf L_tip}-\ref{eq: bf M_tip}) converges to zero as $t$ tends to $+ \infty$.

\section{Simulations}\label{sec: simulations}
The numerical values used in simulation studies are extracted from FW-11. It is a conceptual design of high altitude long endurance commercial aircraft with a wing length $l=16.2\,\mathrm{m}$~\cite{qiao2012effect}. As the original wing of FW-11 is nonhomogeneous, we considered the structural characteristics at the mid-length of the wing, corresponding to the 4th section of the wing presented in~\cite{qiao2012effect}. We set the damping coefficients to $\eta_\phi=\eta_\omega=0.02$ and the coupling offset to $x_c = 2\,\mathrm{m}$. The aerodynamic coefficients provided in~\cite{qiao2012effect} correspond to an altitude of $10,000\,\mathrm{m}$ and a flight speed of $240\,\mathrm{m/s}$. However, they are purely static (i.e., only $\alpha_\omega$ and $\alpha_\phi$ are provided). For simulation purposes, \emph{ad hoc} unsteady aerodynamic coefficients have been set to $\beta_\omega = \gamma_\omega = \alpha_\omega / 3$ and $\beta_\phi = \gamma_\phi = - \alpha_\phi / 2$. Finally, the store characteristics are set to $m_s = 1,000\,\mathrm{kg}$ and $J_s = 500\,\mathrm{kg \cdot m^2}$. With this setup, the conditions in Assumption~\ref{assumption: physical parameters} are satisfied by setting $\varepsilon_1 = 3.189 \times 10^{-4}$ and $\varepsilon_2= 1.195 \times 10^{-3}$.

\begin{figure}[htb]
	\centering
	\includegraphics[width=3in,trim={0 0 0 0}]{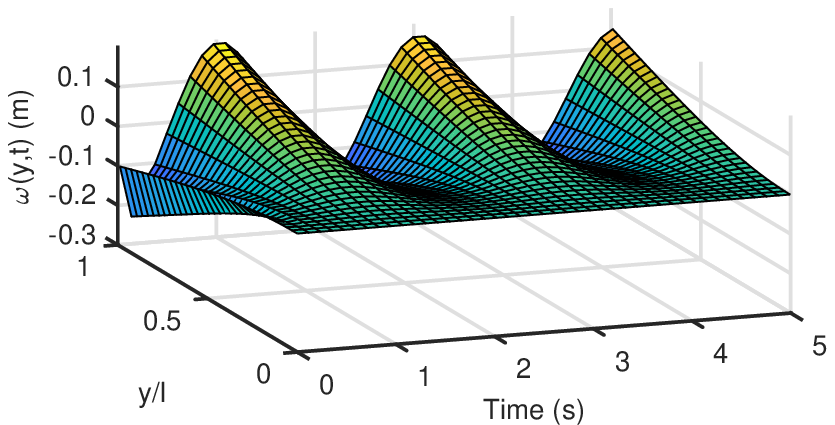}
	\includegraphics[width=3in,trim={0 0 0 0}]{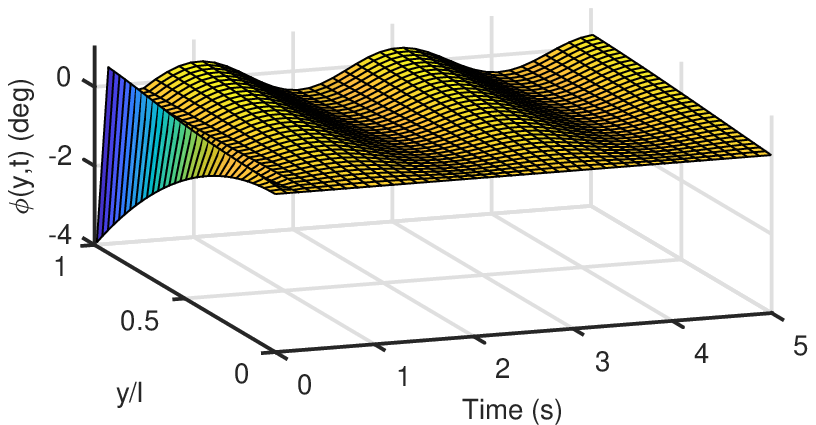}
	\caption{Open-loop response: Bending displacement $\omega(y,t)$; Twisting displacement $\phi(y,t)$.}
	\label{fig: OL}
\end{figure}

\begin{figure}[htb]
	\centering
		\includegraphics[width=3in,trim={0 0 0 0}]{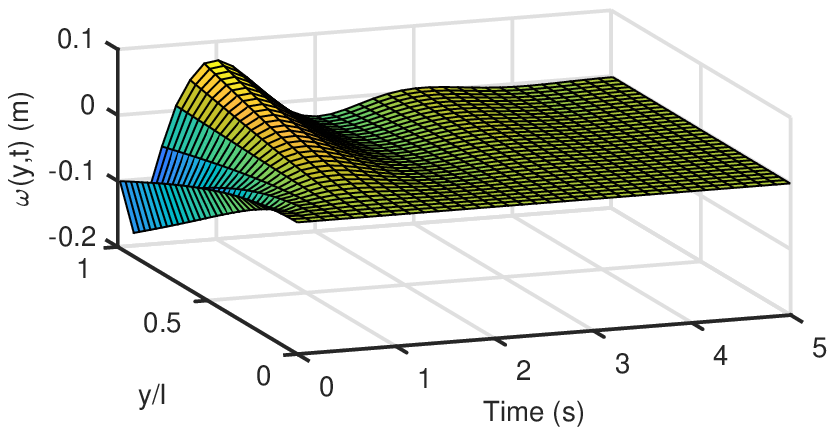}
		\includegraphics[width=3in,trim={0 0 0 0}]{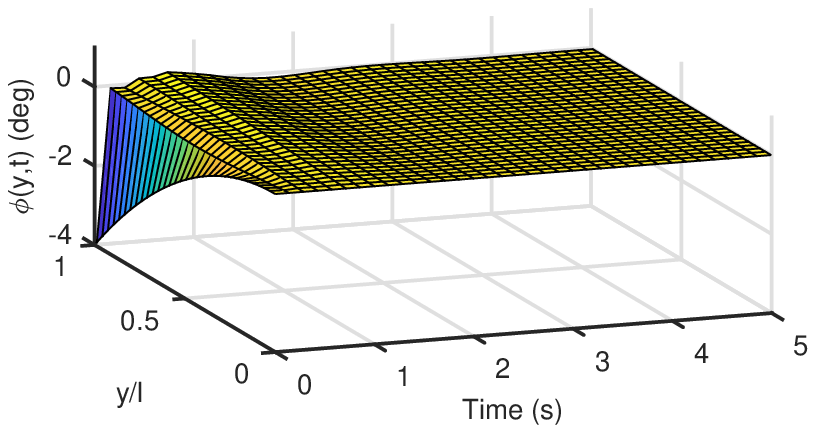}
	\caption{Closed-loop response: Bending displacement $\omega(y,t)$; Twisting displacement $\phi(y,t)$.}
	\label{fig: CL}
\end{figure}

\begin{figure}[htb]
\centering
\includegraphics[width=3.2in]{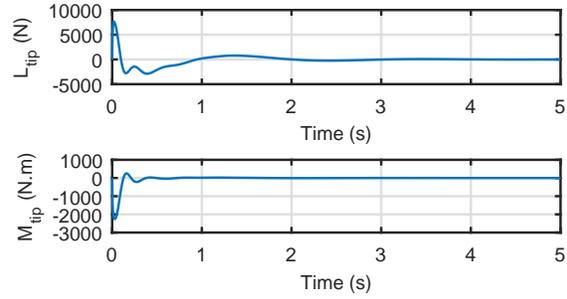}
\caption{Applied control momentum $M_\mathrm{tip}$ and force $F_\mathrm{tip}$}
\label{fig: applied control}
\end{figure}

Numerical simulations are carried out based on the Galerkin method \cite{ciarlet2002finite}. The temporal behavior of the open-loop system is depicted in Fig.~\ref{fig: OL}. It can be seen that both bending and twisting displacements are poorly damped, exhibiting large oscillations. In contrast, by setting the controller gains as $k_1=5,000$ and $k_2=2,500$, the oscillations are damped out rapidly in closed loop as shown in Fig.~\ref{fig: CL}.
The actuation effort at the wing tip is depicted in Fig.~\ref{fig: applied control}.


%

\section{Conclusion}
This paper tackled the boundary control of a flexible wing under unsteady aerodynamic loads and in presence of a store located at the wing tip. The wing is modeled by a distributed parameter system consisting of two coupled partial differential equations describing both bending and twisting displacements. After demonstrating that the problem is well-posed, it is shown by using the Lyapunov method that the proposed boundary control scheme ensures the uniform exponential stability of both bending and twisting displacements. The obtained results rely in certain structural constraints that mainly impose restrictions on the balance between the structural stiffness ant the aerodynamic coefficients. As these constraints may limit the admissible airspeed, evaluating their conservatism and their potential relaxation shall be considered in future work.



\bibliographystyle{plain}        
\bibliography{autosam}           

\begin{thebibliography}{10}

\bibitem{balakrishnan2001subsonic}
A.~V. Balakrishnan.
\newblock Subsonic flutter suppression using self-straining actuators.
\newblock {\em Journal of the Franklin Institute}, 338(2):149--170, 2001.

\bibitem{balakrishnan2003toward}
A.~V. Balakrishnan.
\newblock Toward a mathematical theory of aeroelasticity.
\newblock In {\em IFIP Conference on System Modeling and Optimization}, pages
  1--24. Springer, 2003.

\bibitem{balakrishnan2004spectral}
A.~V. Balakrishnan, M.~A. Shubov, and C.~A. Peterson.
\newblock Spectral analysis of coupled euler-bernoulli and timoshenko beam
  model.
\newblock {\em ZAMM-Journal of Applied Mathematics and Mechanics/Zeitschrift
  f{\"u}r Angewandte Mathematik und Mechanik}, 84(5):291--313, 2004.

\bibitem{beran2004studies}
P.~S. Beran, T.~W. Strganac, K.~Kim, and C.~Nichkawde.
\newblock Studies of store-induced limit-cycle oscillations using a model with
  full system nonlinearities.
\newblock {\em Nonlinear Dynamics}, 37(4):323--339, 2004.

\bibitem{Bhoir2003}
N.~Bhoir and S.~N. Singh.
\newblock Output feedback nonlinear control of an aeroelastic system with
  unsteady aerodynamics.
\newblock {\em Aerospace Science and Technology}, 8(3):195--205, 2004.

\bibitem{Bialy2016}
B.~J. Bialy, I.~Chakraborty, S.~C. Cekic, and W.~E. Dixon.
\newblock Adaptive boundary control of store induced oscillations in a flexible
  aircraft wing.
\newblock {\em Automatica}, 70:230--238, 2016.

\bibitem{Cannarsa2015}
P.~Cannarsa and T.~D'Aprile.
\newblock {\em Introduction to measure theory and functional analysis},
  volume~89.
\newblock Springer, 2015.

\bibitem{chueshov2016nonlinear}
I.~Chueshov, E.~H. Dowell, I.~Lasiecka, and J.~T. Webster.
\newblock Nonlinear elastic plate in a flow of gas: Recent results and
  conjectures.
\newblock {\em Applied Mathematics \& Optimization}, 73(3):475--500, 2016.

\bibitem{chueshov2012generation}
I.~Chueshov and I.~Lasiecka.
\newblock Generation of a semigroup and hidden regularity in nonlinear subsonic
  flow-structure interactions with absorbing boundary conditions.
\newblock {\em Jour. Abstr. Differ. Equ. Appl}, 3:1--27, 2012.

\bibitem{ciarlet2002finite}
P.~G. Ciarlet.
\newblock {\em The finite element method for elliptic problems}.
\newblock SIAM, 2002.

\bibitem{Curtain2009}
R.~Curtain and K.~Morris.
\newblock Transfer functions of distributed parameter systems: A tutorial.
\newblock {\em Automatica}, 45:1101--1116, 2009.

\bibitem{Curtain2012}
R.~F. Curtain and H.~Zwart.
\newblock {\em An introduction to infinite-dimensional linear systems theory},
  volume~21.
\newblock Springer Science \& Business Media, 2012.

\bibitem{Queiroz2012}
M.~S. De~Queiroz, D.~M. Dawson, S.~P. Nagarkatti, and F.~Zhang.
\newblock {\em Lyapunov-based control of mechanical systems}.
\newblock Springer Science \& Business Media, 2012.

\bibitem{Queiroz2002}
M.~S. De~Queiroz and C.~D. Rahn.
\newblock Boundary control of vibration and noise in distributed parameter
  systems: an overview.
\newblock {\em Mechanical Systems and Signal Processing}, 16(1):19--38, 2002.

\bibitem{Guo2002}
B.-Z. Guo.
\newblock Riesz basis approach to the tracking control of a flexible beam with
  a tip rigid body without dissipativity.
\newblock {\em Optimization Methods and Software}, 17(4):655--681, 2002.

\bibitem{Hardy1952}
G.~H. Hardy, J.~E. Littlewood, and G.~P{\'o}lya.
\newblock {\em Inequalities}.
\newblock Cambridge university press, 1952.

\bibitem{He2017}
W.~He and S.~Zhang.
\newblock Control design for nonlinear flexible wings of a robotic aircraft.
\newblock {\em IEEE Transactions on Control Systems Technology},
  25(1):351--357, 2017.

\bibitem{Ko1999}
J.~Ko, T.~W. Strganac, and A.~J. Kurdila.
\newblock Adaptive feedback linearization for the control of a typical wing
  section with structural nonlinearity.
\newblock {\em Nonlinear Dynamics}, 18(3):289--301, 1999.

\bibitem{krstic2008output}
M.~Krstic, B.-Z. Guo, A.~Balogh, and A.~Smyshlyaev.
\newblock Output-feedback stabilization of an unstable wave equation.
\newblock {\em Automatica}, 44(1):63--74, 2008.

\bibitem{Krstic2008}
M.~Krstic and A.~Smyshlyaev.
\newblock {\em Boundary control of PDEs: A course on backstepping designs}.
\newblock SIAM, 2008.

\bibitem{lasiecka2016feedback}
I.~Lasiecka and J.~T. Webster.
\newblock Feedback stabilization of a fluttering panel in an inviscid subsonic
  potential flow.
\newblock {\em SIAM Journal on Mathematical Analysis}, 48(3):1848--1891, 2016.

\bibitem{Leoni2009}
G.~Leoni.
\newblock {\em A first course in Sobolev spaces}, volume 105.
\newblock American Mathematical Society Providence, RI, 2009.

\bibitem{Luo2012}
Z.-H. Luo, B.-Z. Guo, and O.~Morgul.
\newblock {\em Stability and stabilization of infinite dimensional systems with
  applications}.
\newblock Springer Science \& Business Media, 2012.

\bibitem{Mukhopadhyay2000Transonic}
V.~Mukhopadhyay.
\newblock Transonic flutter suppression control law design and wind tunnel test
  results.
\newblock {\em Journal of Guidance, Control, and Dynamics}, 23(5):930--937,
  2000.

\bibitem{Mukhopadhyay2003}
V.~Mukhopadhyay.
\newblock Historical perspective on analysis and control of aeroelastic
  responses.
\newblock {\em Journal of Guidance, Control, and Dynamics}, 26(5):673--684,
  2003.

\bibitem{Paranjape2013}
A.~A. Paranjape, J.~Guan, S.-J. Chung, and M.~Krstic.
\newblock {PDE} boundary control for flexible articulated wings on a robotic
  aircraft.
\newblock {\em IEEE Transactions on Robotics}, 29(3):625--640, 2013.

\bibitem{Pazy2012}
A.~Pazy.
\newblock {\em Semigroups of linear operators and applications to partial
  differential equations}, volume~44.
\newblock Springer Science \& Business Media, 2012.

\bibitem{qiao2012effect}
Yuqing Qiao.
\newblock Effect of wing flexibility on aircraft flight dynamics.
\newblock Master's thesis, Cranfield University, UK, 2012.

\bibitem{Royden1988}
H.~L. Royden and P.~Fitzpatrick.
\newblock {\em Real analysis}, volume 198.
\newblock Macmillan New York, 1988.

\bibitem{Scott2000}
R.~C. Scott, S.~T. Hoadley, C.~D. Wieseman, and M.~H. Durham.
\newblock Benchmark active controls technology model aerodynamic datai.
\newblock {\em Journal of Guidance, Control, and Dynamics}, 23(5):914--921,
  2000.

\bibitem{Shearer2007}
C.~M. Shearer and C.~E. Cesnik.
\newblock Nonlinear flight dynamics of very flexible aircraft.
\newblock {\em Journal of Aircraft}, 44(5):1528--1545, 2007.

\bibitem{stanewsky2001adaptive}
E~Stanewsky.
\newblock Adaptive wing and flow control technology.
\newblock {\em Progress in Aerospace Sciences}, 37(7):583--667, 2001.

\bibitem{Su2010}
W.~Su and C.~E. Cesnik.
\newblock Nonlinear aeroelasticity of a very flexible blended-wing-body
  aircraft.
\newblock {\em Journal of Aircraft}, 47(5):1539--1553, 2010.

\bibitem{zhang2011spectrum}
G.-D. Zhang and B.-Z. Guo.
\newblock On the spectrum of euler-bernoulli beam equation with kelvin-voigt
  damping.
\newblock {\em Journal of Mathematical Analysis and Applications},
  374(1):210--229, 2011.

\bibitem{Zhang2005}
X.~Zhang, W.~Xu, S.~S. Nair, and V.~Chellaboina.
\newblock {PDE} modeling and control of a flexible two-link manipulator.
\newblock {\em IEEE Transactions on Control Systems Technology},
  13(2):301--312, 2005.

\bibitem{Ziabar2010}
M.~Y. Ziabari and B.~Ghadiri.
\newblock Vibration analysis of elastic uniform cantilever rotor blades in
  unsteady aerodynamics modeling.
\newblock {\em Journal of Aircraft}, 47(4):1430--1434, 2010.

\end{thebibliography}



\end{document}